\documentclass[reqno,11pt]{amsart}
\usepackage{amssymb}
\usepackage{latexsym,array}
\usepackage{amsfonts}
\newcommand{\C}{\mathbb {C}}
\newcommand{\D}{{\mathbb D}}

\newcommand{\sbm}[1]{\left[\begin{smallmatrix} #1
                \end{smallmatrix}\right]}
\newcommand{\bl}{\boldsymbol{e_\ell}}
\newcommand{\br}{\boldsymbol{e_r}}

\newcommand{\bH}{\mathbb{H}}

\newtheorem{Pa}{Paper}[section]
\newtheorem{theorem}[Pa]{{\bf Theorem}}

\newtheorem{definition}[Pa]{{\bf Definition}}

\newtheorem{remark}[Pa]{{\bf Remark}}

\author[V. Bolotnikov]{Vladimir Bolotnikov}
\address{Department of Mathematics\\
The College of William and Mary \\  
Williamsburg, VA 23187-8795\\
USA}
\email{vladi@math.wm.edu}
\title[Pick matrices]{Pick matricies and 
quaternionic power series}
\begin{document}
\begin{abstract}
It is well known that a non-constant complex-valued function $f$ defined 
on the open unit disk $\mathbb D$ is an analytic self-mapping of $\D$
if and only if Pick matrices $\left[
(1-f(z_i)\overline{f(z_j)})/(1-z_i\overline{z}_j)\right]_{i,j=1}^n$ are positive semidefinite
for all choices of finitely many points $z_i\in\D$. A stronger version of the ``if" part was established 
by A. Hindmarsh \cite{hind}: if all $3\times 3$ Pick matrices are positive semidefinite, then
$f$ is an analytic self-mapping of $\D$. In this paper, we extend this result to the non-commutative setting of 
power series over quaternions.
\end{abstract}
\maketitle

\maketitle 
\section{Introduction}
\setcounter{equation}{0}

The Schur class ${\mathcal S}$ of all analytic  complex-valued functions mapping the 
unit disk $\D$ into its closure has played a prominent role in function theory and its
applications beginning with the work of I. Schur \cite{Schur}.
Among several alternative characterizations of the Schur class is
one in terms of positive kernels:  
{\em the function $f \colon {\mathbb D} \to {\mathbb
C}$ is in the class ${\mathcal S}$ if and only if the
associated kernel
$K_f(z,\zeta)=\frac{1-f(z)\overline{f(\zeta)}}{1-z\bar{\zeta}}$
is positive on $\D\times\D$ or equivalently, if and only if the 
{\em Pick matrix}
\begin{equation}
P_f(z_1,\ldots,z_n)=
\left[\frac{1-f(z_i)\overline{f(z_j)}}{1-z_i\overline{z}_j}\right]_{i,j=1}^n
\label{1.1}
\end{equation}
is positive semidefinite for any choice of finitely many points 
$z_1,\ldots,z_n\in\D$}. The ``only if" part is the classical result
of Pick and Nevanlinna \cite{pick, nevan1}. For the ``if" part, 
let us observe that positivity of $1\times 1$ matrices
$P_f(z)$ already guarantees $|f(z)|\le 1$ ($z\in\D$). Thus, larger Pick matrices
are needed in the ``if" direction to guarantee the analyticity of $f$.
The latter can be done via constructing the coisometric de Branges-Rovnyak realization 
\cite{dbr2} for $f$ or using a more recent lurking isometry argument 
\cite{Ball-Winnipeg}. A remarkable fact established by Hindmarsh \cite{hind}
(see also \cite{fh}) is that analyticity is implied by positivity of all $3\times 3$ Pick 
matrices. The objective of this paper is to extend this result to regular
functions in quaternionic variable (Theorem \ref{T:1.3} below).

\smallskip

Let $\mathbb H$ denote the skew field of real quaternions
$\alpha=x_0+{\bf i}x_1+{\bf j}x_2+{\bf k}x_3$  where $x_\ell\in\mathbb R$ and
${\bf i}, {\bf j}, {\bf k}$ are imaginary units such that
${\bf i}^2={\bf j}^2={\bf k}^2={\bf ijk}=-1$. The real and the imaginary
parts, the conjugate and the absolute value of a quaternion $\alpha$   are defined by
${\rm Re}(\alpha)=x_0$, ${\rm Im}(\alpha)=ix_1+jx_2+kx_3$, $\bar \alpha={\rm 
Re}(\alpha)-{\rm Im}(\alpha)$ and 
$|\alpha|^2=\alpha\bar{\alpha}=|{\rm Re}(\alpha)|^2+|{\rm Im}(\alpha)|^2$, 
respectively. By $\mathbb B=\{\alpha\in\mathbb H: \; |\alpha|<1\}$ we denote the unit ball
in $\mathbb H$.

\smallskip

Since multiplication in $\mathbb H$ is not commutative, function theory 
over quaternions is somewhat different from that over the field $\C$.
There have been several notions of regularity (or analyticity) for $\mathbb H$-valued 
functions, most notable of which are due to Moisil \cite{moisil}, Fueter \cite{fuet1, fuet2}, 
and Cullen \cite{cullen}. More recently, upon refining and developing Cullen's approach, 
Gentili and Struppa introduced in \cite{genstr} the notion of regularity which,
being restricted to functions on a quaternionic ball around the origin, turns out to be the 
feature of power series with quaternionic coefficients on one side; we refer to the recent book
\cite{gss} for a detailed exposition of the subject. Here we  accept the 
following definition of regularity on the quaternionic unit ball.
\begin{definition}
{\rm A function $f: \, \mathbb B\to \mathbb H$ is called {\em left-regular} on $\mathbb B$ 
if it admits the power series expansion with quaternionic coefficients on the right which 
converges absolutely on $\mathbb B$:
\begin{equation}
f(z)=\sum_{k=0}^\infty z^kf_k\quad\mbox{with $\; \; f_k\in\mathbb H\; $ such that $\; \;
\overline{\displaystyle\lim_{k\to\infty}}\sqrt[k]{|f_k|}\le 1$}.
\label{2.2}   
\end{equation}
If in addition, $f(\alpha):={\displaystyle\sum_{k=0}^\infty \alpha^kf_k}\in \overline{\mathbb B}$ 
for all  $\alpha\in\mathbb B$, we say 
that $f$ belongs to the {\em left Schur class} $\mathcal {QS}_{_L}$. 
Right regular functions and the {\em right Schur class} $\mathcal {QS}_{_R}$ can be defined 
similarly.}
\label{D:1.1}
\end{definition}
Quaternionic Schur classes have become an object of intensive study quite recently. 
A number of related results (e.g., M\"obius transformations, Schwarz Lemma, Bohr's inequality)
are presented in \cite[Chapter 9]{gss}. Among other results, we mention realizations for slice 
regular functions \cite{acs1}, Schwarz-Pick Lemma \cite{bist}, Blaschke products \cite{acs2},
Nevanlinna-Pick interpolation \cite{abcs}.

\section{Pick matrices and Hindmarsh's theorem}                        
\setcounter{equation}{0}

A straightforward entry-wise verification shows that the complex matrix \eqref{1.1} 
satisfies the Stein equality 
\begin{equation} 
P_f(z_1,\ldots,z_n)-TP_f(z_1,\ldots,z_n)T^*=EE^*-NN^*,
\label{2.3}
\end{equation}
where 
\begin{equation}
T=\begin{bmatrix}z_1 & & 0 \\ & \ddots & \\ 0 && z_n\end{bmatrix},\quad
E=\begin{bmatrix}1 \\ \vdots \\ 1\end{bmatrix},\quad N=\begin{bmatrix}f(z_1) \\ \vdots \\ 
f(z_n)\end{bmatrix}.
\label{2.3a}   
\end{equation}
Since $|z_i|<1$, the latter matrix is the {\em unique} matrix subject to identity \eqref{2.3}. 
In case $z_i\in\mathbb B$ and $f(z_i)\in\mathbb H$,
the Stein equation \eqref{2.3} still has a unique solution $P_f(z_1,\ldots,z_n)$ 
(still called the Pick matrix of $f$). Solving this equation gives the 
explicit formula for the entries of 
$P_f(z_1,\ldots,z_n)$ in terms of series
\begin{equation} 
P_f(z_1,\ldots,z_n)=
\left[\sum_{k=0}^\infty z_i^k(1-f(z_i)\overline{f(z_j)})\overline{z}_j^k\right]_{i,j=1}^n
\label{1.4}
\end{equation}
which converge due to the following estimate:
$$
\left|\sum_{k=0}^\infty z_i^k(1-f(z_i)\overline{f(z_j)})\overline{z}_j^k\right|
\le 2\sum_{k=0}^\infty | z_i|^k |z_j|^k=\frac{2}{1-| z_i||z_j|}.
$$
According to a result from \cite{abcs}, for any function $f\in\mathcal {QS}_{_L}$, the 
associated Pick matrix \eqref{1.4}
is positive semidefinite for any choice of finitely many points
$z_1,\ldots,z_n\in\mathbb B$.  The notions of adjoint matrices, of 
Hermitian matrices and positive semidefinite matrices over  $\mathbb H$ are similar
to those  over $\mathbb C$ (we refer to a very nice survey \cite{fuzhen} on this
subject).

\smallskip

The following quaternionic analog of the Hindmarsh 
theorem \cite{hind} is the main result of the present paper.
\begin{theorem}
Let $f: \, \mathbb B\to \mathbb H$ be given and let us assume that $3\times 3$
Pick matrices $P_f(z_1,z_2,z_3)$ are positive semidefinite for all 
$(z_1,z_2,z_3)\in\mathbb B^3$. Then $f$ belongs to $\mathcal {QS}_{_L}$.
\label{T:1.3}
\end{theorem}
Before starting the proof we make several observations. 
\begin{remark}
For any $z_1,z_2\in\C$, the quaternion $z_1{\bf j}z_2$ belongs to 
$\C{\bf j}$.
\label{R:2.6}
\end{remark}
\noindent
The statement follows from the multiplication table for imaginary units in $\mathbb H$.
We also remark that any quaternion $\alpha=x_0+{\bf i}x_1+{\bf j}x_2+{\bf k}x_3$ admits 
a unique representation $\alpha=\alpha_1+\alpha_2{\bf
j}$ with $\alpha_1,\alpha_2\in\C$. 
%(in fact, $\alpha_1=x_0+{\bf i}x_1$ and $\alpha_2=x_2+{\bf i}x_3$).
Consequently, any quaternionic matrix $A$ admits a 
unique representation $A=A_1+A_2{\bf j}$ with complex matrices $A_1$ and $A_2$.
\begin{remark}
The matrix $A=A_1+A_2{\bf j}\in\mathbb H^{n\times n}$ ($A_1,A_2\in\C^{n\times n}$)
is positive semidefinite if and only if the complex matrix 
$\begin{bmatrix}A_1 & A_2 \\ -\overline{A}_2 & \overline{A}_1\end{bmatrix}$ is positive semidefinite
(see \cite{fuzhen}).
\label{R:2.2}
\end{remark}
{\rm Two quaternions $\alpha$ and $\beta$ are called {\em equivalent} (conjugate to each other)
if $\alpha=h^{-1}\beta h$ for some nonzero $h\in\mathbb H$.} It follows (see e.g., \cite{fuzhen}) 
that
\begin{equation}
\alpha\sim\beta\quad\mbox{if and only if}\quad {\rm Re}(\alpha) ={\rm Re}(\beta) \; \;
\mbox{and} \; \; |\alpha|=|\beta|.
\label{2.0}
\end{equation}
Therefore, the conjugacy class of a given $\alpha\in\mathbb H$ form a $2$-sphere (of radius
$|{\rm Im}(\alpha)|$ around ${\rm Re}(\alpha)$). 
\begin{remark}
If $\alpha,\beta,\gamma$ are three distinct equivalent quaternions, then
\begin{equation}
\gamma^k=(\gamma-\beta)(\alpha-\beta)^{-1}\alpha^k+(\alpha-\gamma)(\alpha-\beta)^{-1}\beta^k
\quad\mbox{for all}\quad k=0,1,\ldots
\label{2.1}
\end{equation}
{\rm Indeed, since $\alpha\sim\beta$, it follows from \eqref{2.0} that
$(\alpha-\beta)\alpha(\alpha-\beta)^{-1}=\overline{\beta}$ and
subsequently,
\begin{equation}
(\alpha-\gamma)a^k(\alpha-\gamma)^{-1}=\overline{\beta}^k=(\gamma-\beta)^{-1}\gamma^k(\gamma-\beta),
\label{2.1a}  
\end{equation}
where the first equality is clear and the second is a virtue of the first
since $\beta$ and $\gamma$ are also equivalent. Similarly,
\begin{equation}
(\alpha-\beta)^{-1}\beta^k(\alpha-\beta)=\overline{a}^k=(\alpha-\gamma)^{-1}\gamma^k(\alpha-\gamma)
\label{2.1b}  
\end{equation}
for all integers $k\ge 0$. Then we get \eqref{2.1} from \eqref{2.1a} and \eqref{2.1b} as follows:
\begin{align*}
&(\gamma-\beta)(\alpha-\beta)^{-1}\alpha^k+(\alpha-\gamma)(a-\beta)^{-1}\beta^k\\
&=\gamma^k(\gamma-\beta)(\alpha-\beta)^{-1}+\gamma^k(\alpha-\gamma)(\alpha-\beta)^{-1}=\gamma^k.
\end{align*}}   
\label{R:2.0}
\end{remark}
It turns out that the values of a regular function $f$ at two points from the same conjugacy class
uniquely determine $f$ at any point from this class; the formula \eqref{2.4} below was established 
in \cite{genstr} in a more general setting. 
\begin{remark}
{\rm Let $f$ be left-regular on $\mathbb B$ and let $\alpha,\beta,\gamma\in\mathbb B$ be distinct
equivalent points. Then
\begin{equation}
f(\gamma)=(\gamma-\beta)(\alpha-\beta)^{-1}f(\alpha)+(\alpha-\gamma)(\alpha-\beta)^{-1}f(\beta). 
\label{2.4}
\end{equation}  
Indeed, equality \eqref{2.1} verifies formula \eqref{2.4}
for monomials $f(z)=z^k$. Extending it by right linearity to power
series \eqref{2.2} completes the proof.}
\label{R:rep}
\end{remark}
\noindent
{\bf Proof of Theorem \ref{T:1.3}:} We first observe that for $n=1$, the formula 
\eqref{1.4} amounts to 
$$
P_f(z_1)=\sum_{k=0}^\infty 
z_1^k(1-|f(z_i)|^2)\overline{z}_1^k=\frac{1-|f(z_1)|^2}{1-|z_1|^2}
$$
for each $z_1\in\mathbb B$. Therefore, condition $P_f(z_1)\ge 0$ implies
$|f(z_1)|\le 1$. 

\smallskip

It remains to show that $f$ is left-regular. Toward this end, we first
show that there exist complex Schur-class functions $g$ and $h$ such that
\begin{equation}
f(\zeta)=g(\zeta)+h(\zeta){\bf j}\quad\mbox{for all}\quad \zeta\in\D.
\label{op}
\end{equation}
Indeed, for each fixed {\em complex} point $\zeta\in\C\cap\mathbb B=\D$, the quaternion
$f(\zeta)\in\mathbb H$ admits a (unique) representation \eqref{op} with 
$g(\zeta)\in\C$ and  $h(\zeta)\in\C$. For any two points 
$\zeta_1,\zeta_2\in\D$ and any $k\ge 0$, we then compute
\begin{align*}
\zeta_1^k(1-f(\zeta_1)\overline{f(\zeta_2)})\overline{\zeta}_2^k
=&\zeta_1^k\overline{\zeta}_2^k-\zeta_1^k\left[g(\zeta_1)+h(\zeta_1){\bf j})
(\overline{g(\zeta_2)}-{\bf j}\overline{h(\zeta_2)}\right]\overline{\zeta}_2^k\\
=&\zeta_1^k\left[1-g(\zeta_1)\overline{g(\zeta_2)}- 
h(\zeta_1)\overline{h(\zeta_2)}\right]\overline{\zeta}_2^k\\
&+\zeta_1^k\left[g(\zeta_1){\bf j}\overline{h(\zeta_2)}-h(\zeta_1){\bf j}\overline{g(\zeta_2)}\right]
\overline{\zeta}_2^k.
\end{align*}
Summing up the latter equalities over all $k\ge 0$ gives
\begin{align}
\sum_{k=0}^\infty \zeta_1^k(1-f(\zeta_1)\overline{f(\zeta_2)})\overline{\zeta}_2^k=&
\frac{1-g(\zeta_1)\overline{g(\zeta_2)}-h(\zeta_1)\overline{h(\zeta_2)}}{1-\zeta_1\overline{\zeta}_2}
\label{2.6}\\
&+\sum_{k=0}^\infty \zeta_1^k\left[g(\zeta_1){\bf j}\overline{h(\zeta_2)}-h(\zeta_1){\bf 
j}\overline{g(\zeta_2)}\right]\overline{\zeta}_2^k.\notag
\end{align}
The first term on the right is complex whereas the second term belongs to $\C{\bf j}$ by Remark 
\ref{R:2.6}. Let us consider the Pick matrix $P_f(\zeta_1,\zeta_2,\zeta_3)$ based on arbitrary
points $\zeta_1,\zeta_2,\zeta_3\in\C$.  According to \eqref{1.4} and \eqref{2.6},
$$
P_f(\zeta_1,\zeta_2,\zeta_3)=P_{f,1}(\zeta_1,\zeta_2,\zeta_3)+P_{f,2}(\zeta_1,\zeta_2,\zeta_3)
$$
where
\begin{equation}
P_{f,1}(\zeta_1,\zeta_2,\zeta_3)=
\left[\frac{1-g(\zeta_i)\overline{g(\zeta_j)}-h(\zeta_i)\overline{h(\zeta_j)}}
{1-\zeta_i\overline{\zeta}_j}\right]_{i,j=1}^3\in\C^{3\times 3}
\label{2.7}
\end{equation}
and where $P_{f,2}(\zeta_1,\zeta_2,\zeta_3)$ is a matrix from $\C^{3\times 3}{\bf j}$.
By the assumption in Theorem \ref{T:1.3}, the matrix $P_f(\zeta_1,\zeta_2,\zeta_3)$ is positive 
semidefinite. Then, the complex matrix $P_{f,1}(\zeta_1,\zeta_2,\zeta_3)$ is 
positive semidefinite by 
Remark \ref{R:2.2}. The matrix \eqref{2.7} can be written as 
\begin{equation}
P_{f,1}(\zeta_1,\zeta_2,\zeta_3)=\Lambda-G\Lambda G^*-H\Lambda H^*\ge 0
\label{2.8}
\end{equation}
where 
$$
\Lambda=\left[\frac{1}
{1-\zeta_i\overline{\zeta}_j}\right]_{i,j=1}^3,\quad
G=\sbm{g(\zeta_1) & 0 & 0 \\ 0 & g(\zeta_2) & 0 \\ 0 & 0 & g(\zeta_3)},\quad 
H=\sbm{h(\zeta_1) & 0 & 0 \\ 0 & h(\zeta_2) & 0 \\ 0 & 0 & h(\zeta_3)},
$$
and it is well known that the matrix $\Lambda$ is positive semidefinite. 
Then we conclude from \eqref{2.8} that 
$\Lambda-G\Lambda G^*\ge 0$ and $\Lambda-H\Lambda H^*\ge 0$, i.e., that the $3\times 3$ matrices 
$$
\left[\frac{1-g(\zeta_i)\overline{g(\zeta_j)}}
{1-\zeta_i\overline{\zeta}_j}\right]_{i,j=1}^3\quad\mbox{and}\quad
\left[\frac{1-h(\zeta_i)\overline{h(\zeta_j)}}
{1-\zeta_i\overline{\zeta}_j}\right]_{i,j=1}^3
$$
are positive semidefinite for all choices of $\zeta_1,\zeta_2,\zeta_3\in\D$. Then 
it follows from 
the complex Hindmarsh theorem that the functions $g,h: \, \D\to\C$ are complex-analytic
and belong to the classical Schur class $\mathcal S$. Substituting their power series expansions 
$g(\zeta)={\displaystyle\sum_{k=0}^\infty\zeta^k g_k}$ and $h(\zeta)={\displaystyle\sum_{k=0}^\infty\zeta^k 
h_k}$ into \eqref{op} leads us to the power series expansion for $f$ on $\D$:
\begin{equation}
f(\zeta)=g(\zeta)+ h(\zeta){\bf j}=\sum_{k=0}^\infty \zeta^k f_k\quad\mbox{with}\quad f_k=g_k+h_k{\bf j}.
\label{2.8a}
\end{equation}
We now extend the latter power series to the whole $\mathbb B$ by simply letting
\begin{equation}
F(z)=\sum_{k=0}^\infty z^k f_k\quad (z\in\mathbb B).
\label{2.9}
\end{equation}
The resulting power series converges absolutely on $\mathbb B$ (as the sum of two converging series 
$g(z)$ and $h(z){\bf j}$) and agrees with $f$ on $\D$. We next show that 
$F$  agrees  with $f$ throughout $\mathbb B$.

\smallskip

Let $\gamma$ be any point in $\mathbb B\setminus\D$. The points $\alpha:={\rm Re}(\gamma)+|{\rm 
Re}(\gamma)|{\bf i}$ and $\overline{\alpha}$ belong to $\D$ and are equivalent to $\gamma$.
Observe that
\begin{equation}
\begin{array}{ll}
(\gamma-\overline{\alpha})(\alpha-\overline{\alpha})^{-1}&=(\gamma-\overline{\gamma})^{-1}(\gamma-\overline{\alpha}),
\\ (\alpha-\gamma)(\alpha-\overline{\alpha})^{-1}&=(\gamma-\overline{\gamma})^{-1}(\gamma-\alpha).
\end{array}
\label{2.9a}
\end{equation}
Since $F$ is left-regular by construction, we apply formula \eqref{2.4} to get
\begin{align}
F(\gamma)=&(\gamma-\overline{\alpha})(a-\overline{\alpha})^{-1}F(\alpha)+(\alpha-\gamma)(\alpha-\overline{\alpha})^{-1}
F(\overline{\alpha})\notag\\
=&(\gamma-\overline{\gamma})^{-1}(\gamma-\overline{\alpha})f(\alpha)+(\gamma-\overline{\gamma})^{-1}(\gamma-\alpha)f(\overline{\alpha}),
\label{2.10}
\end{align} 
where the second equality follows due to \eqref{2.9a} and 
since $F$ agrees with $f$ on $\D$. On the other hand, we know that the Pick matrix 
$P_f(\alpha,\overline{\alpha},c)$ is 
positive semidefinite and 
satisfies the Stein identity \eqref{2.3}:
\begin{equation}
P_f(\alpha,\overline{\alpha},\gamma)-TP_f(\alpha,\overline{\alpha},\gamma)T^*=EE^*-NN^*
\label{2.11}
\end{equation}
where 
$$
T=\begin{bmatrix} \alpha & 0 & 0 \\ 0 & \overline{\alpha} & 0 \\ 0 & 0 & \gamma\end{bmatrix},\quad
E=\begin{bmatrix} 1 \\ 1\\ 1\end{bmatrix},\quad N=\begin{bmatrix} f(\alpha) \\ f(\overline{a})  \\ 
f(\gamma)\end{bmatrix}.
$$
Let us introduce the row-vector $V=\begin{bmatrix} \gamma-\overline{a} & &\gamma-\alpha && 
\overline{\gamma}-\gamma\end{bmatrix}$. Since ${\rm Re}(\alpha)={\rm Re}(\gamma)$ and 
$|\alpha|=|\gamma$ by definition of $\alpha$, we have
\begin{align*}
VT&=\begin{bmatrix} (\gamma-\overline{\alpha})\alpha & (\gamma-\alpha)\overline{\alpha} & 
(\overline{\gamma}-\gamma)\gamma\end{bmatrix}=\gamma V,\\
VE&=\gamma-\overline{\alpha}+\gamma-\alpha+\overline{\gamma}-\gamma=2{\rm Re}(\gamma)-2{\rm Re}(\alpha)=0.
\end{align*}
Multiplying both parts of \eqref{2.11} by $V$ on the left and by $V^*$ on the right we get,
on account of two last equalities,
$$
VP_f(\alpha,\overline{\alpha},\gamma)V^*-\gamma 
VP_f(\alpha,\overline{\alpha},\gamma)V^*\overline{\gamma}=-VNN^*V^*.
$$
Since $P_f(\alpha,\overline{\alpha},\gamma)$ is positive semidefinite, we have 
$VP_f(\alpha,\overline{\alpha},\gamma)V^*\ge 0$
and hence we can write the last equality as 
$$
(1-|c|^2)VP_f(\alpha,\overline{\alpha},\gamma)V^*=-|VN|^2.
$$
The latter may occur only if $VP_f(\alpha,\overline{\alpha},\gamma)V^*=VN=0$. Thus,
$$
VN=(c-\overline{\alpha})f(a)+(c-\alpha)f(\overline{\alpha})+(\overline{\gamma}-\gamma)f(\gamma)=0,
$$
which implies
\begin{equation}
f(\gamma)=(\gamma-\overline{\gamma})^{-1}(\gamma-\overline{\alpha})f(\alpha)+(\gamma-\overline{\gamma})^{-1}
(\gamma-\alpha)f(\overline{\alpha}).
\label{2.12}
\end{equation}  
Comparing \eqref{2.10} and \eqref{2.12} we conclude that $F(\gamma)=f(\gamma)$. Since $\gamma$ was chosen 
arbitrarily
in $\mathbb B\setminus\D$, it follows that $F=f$ on $\mathbb B$. Since $F$ is left-regular on 
$\mathbb B$ by construction \eqref{2.9}, it follows that $f$ is left-regular on $\mathbb B$ as well.
\qed

\section{Schur-class of quaternionic formal power series}
\setcounter{equation}{0}   

It turns out that the quaternionic Schur class can be defined without distinguishing the 
left and the right settings. Let us consider formal power series in one
formal variable $z$ which commutes with quaternionic coefficients (which in turn, satisfy 
the same growth condition as in \eqref{2.2}):
\begin{equation}
g(z)=\sum_{k=0}^\infty z^kg_k=\sum_{k=0}^\infty g_k z^k
\quad\mbox{with $\; \; g_k\in\mathbb H\; $ such that $\; \;
\overline{\displaystyle\lim_{k\to\infty}}\sqrt[k]{|g_k|}\le 1$}.
\label{3.1}
\end{equation}
For each $g\in\bH[[z]]$ as in \eqref{3.1}, we define its {\em conjugate} by
$\; g^\sharp(z)={\displaystyle\sum_{k=0}^\infty z^k \overline{g}_k}$.
The anti-linear involution $g\mapsto g^\sharp$ can be viewed as an extension of the 
quaternionic conjugation $\alpha\mapsto \overline{\alpha}$ from $\bH$ to $\bH[[z]]$.
We next define $g^{\bl}(\alpha)$ and $g^{\br}(\alpha)$ (left and right evaluations of $g$ at $\alpha$) by
\begin{equation}
g^{\bl}(\alpha)=\sum_{k=0}^\infty\alpha^k g_k,\quad
g^{\br}(\alpha)=\sum_{k=0}^\infty g_k\alpha^k,\quad\mbox{if}\quad g(z)=\sum_{k=0}^n z^k g_k.
\label{3.2}
\end{equation}
Observe that condition $\;\overline{\displaystyle\lim_{k\to\infty}}\sqrt[k]{|g_k|}\le 1$ imposed on the 
coefficients  guarantees the absolute convergence of the series in \eqref{3.2} for all $\alpha\in\bH$.
Since multiplication in $\bH$ is not commutative, left and right evaluations produce different results;
however, equality $g^{\br}(\alpha)=\overline{g^{\sharp\bl}(\overline{\alpha})}$ holds for any $\alpha\in\bH$
as can be seen from \eqref{3.2} and the definition of $g^\sharp$.

\smallskip

In accordance with Definition \eqref{D:1.1} we define the left and the right Schur classes
$\mathcal {QS}_{_L}$ and $\mathcal {QS}_{_R}$ as the sets of
power series $f\in\bH[[z]]$ such that $|g^{\bl}(\alpha)|\le 1$ (respectively, $|g^{\br}(\alpha)|\le 1$) for all
$\alpha\in\mathbb B$. But as was shown in \cite{abcs} (in slightly different terms), the classes $\mathcal
{QS}_{_L}$ and $\mathcal {QS}_{_R}$ coincide. We now recall several results from \cite{abcs} in terms of 
the present setting. With a power series $g\in\bH[[z]]$ as in \eqref{3.1}, we associate lower triangular Toepliz 
matrices
\begin{equation}
T_n(g)=\begin{bmatrix}g_{0} & 0 & \ldots & 0
\\ g_{1}& g_{0} & \cdot & \cdot \\ \cdot& \cdot & \cdot & 0
\\ g_{n-1}&  \ldots & g_{1} & g_{0}\end{bmatrix}\quad\mbox{for}\quad n=1,2\ldots.
\label{3.3}   
\end{equation}
\begin{theorem}
Let $g\in \bH[[z]]$ be as in \eqref{3.1}. The following are equivalent:
\begin{enumerate}
\item $|g^{\bl}(\alpha)|\le 1$ for all $\alpha\in\mathbb B$.
\item $|g^{\br}(\alpha)|\le 1$ for all $\alpha\in\mathbb B$.
\item The matrix $T_n(g)$ is contractive for all $n\ge 1$.
\end{enumerate}
\label{T:3.1}
\end{theorem}
We thus may talk about the Schur class $\mathcal {QS}\subset\bH[[z]]$ of formal power series $g$
such that the matrix $T_n(g)$ is contractive for all $n\ge 0$. In the latter power series setting, Theorem 
\ref{T:1.3} can be formulated as follows: {\em if the function $f: \, \mathbb B\to \mathbb H$ is such that 
$3\times 3$ Pick matrices $P_f(z_1,z_2,z_3)$ are positive semidefinite for all
$(z_1,z_2,z_3)\in\mathbb B^3$, then there is (a unique) $g\in\mathcal {QS}$ such that 
$f(\alpha)=g^{\bl}(\alpha)$ for all $\alpha\in\mathbb B$.} The ``right" version of this theorem 
is based on dual Pick matrices
$$
\widetilde{P}_f(z_1,\ldots,z_n)=
\left[\sum_{k=0}^\infty \overline{z}_i^k(1-\overline{f(z_i)}f(z_j))z_j^k\right]_{i,j=1}^n.
$$
\begin{theorem}
Let $f: \, \mathbb B\to \mathbb H$ be given and let us assume that $3\times 3$
dual Pick matrices $\widetilde{P}_f(z_1,z_2,z_3)$ are positive semidefinite for all
$(z_1,z_2,z_3)\in\mathbb B^3$. Then there is (a unique) $g\in\mathcal {QS}$ such that
$f(\alpha)=g^{\br}(\alpha)$ for all $\alpha\in\mathbb B$.
\label{T:3.2}
\end{theorem}
The proof is immediate: by Theorem \ref{T:1.3}, there is an $h\in\mathcal {QS}$ such that 
$\overline{f(\overline{\alpha})}=h^{\bl}(\alpha)$ for all $\alpha\in\mathbb B$. Therefore,
$f(\alpha)=\overline{h^{\bl}(\overline{\alpha})}$ and it remains to choose $g=h^\sharp$ 
which belongs to $\mathcal {QS}$ by Theorem \ref{T:3.1}.

\smallskip

In the proof of Theorem \ref{T:1.3}, we actually showed that for any $g\in\mathcal {QS}$,
there exist (unique) Schur-class functions $s,h:\D\to\overline{\D}$ so that 
\begin{equation}
g(\zeta)=s(\zeta)+h(\zeta){\bf j}\quad \mbox{for all}\quad \zeta\in\D,
\label{3.4}
\end{equation}
and the latter equality determines $g$ uniquely in the whole $\mathbb B$. The last question 
we address here is how to characterize the pairs $(s,h)$ of complex Schur functions producing
via formula \eqref{3.4} a quaternionic Schur-class power series. 
\begin{theorem}
Let $s$ and $h$ be Schur-class functions. The function  $g$ be given by \eqref{3.4} belongs to 
$\mathcal {QS}$ if and only if the following matrix is positive semidefinite
\begin{equation}
\begin{bmatrix}{\bf I}_{n}-T_n(s)T_n(s)^*-T_n(h)T_n(h)^* & T_n(s)T_n(h)^\top-T_n(h) T_n(s)^\top\\ \\
\overline{T_n(h)}T_n(s)^*-\overline{T_n(s)}T_n(h)^*& 
{\bf I}_{n}-\overline{T_n(s)}T_n(s)^\top-\overline{T_n(h)}T_n(h)^\top
\end{bmatrix}\ge 0
\label{3.5}   
\end{equation}
for all $n\ge 0$ where ${\bf I}_{n}$ stands for the $n\times n$ identity matrix and $T_n$ is defined
via formula \eqref{3.3}.
\label{T:3.3}
\end{theorem}
{\bf Proof:} By Theorem \ref{T:3.1}, $g$ belongs to $\mathcal {QS}$ if and only if 
${\bf I}_{n}-T_n(g)T_n(g)^*$ is positive semidefinite for all $n\ge 1$. It follows from \eqref{3.4} that 
\begin{align*}
T_n(g)T_n(g)^*&=(T_n(s)+T_n(h){\bf j})(T_n(s)^*-{\bf j}T_n(h)^*)\\
&=T_n(s)T_n(s)^*+T_n(h){\bf j}T_n(s)^*-T_n(s){\bf j}T_n(h)^*+T_n(h)T_n(h)^*\\
&=T_n(s)T_n(s)^*+T_n(h)T_n(h)^*\\
&\quad +(T_n(h)T_n(s)^\top-T_n(s)T_n(h)^\top){\bf j}.
\end{align*}
Therefore, ${\bf I}_{n}-T_n(g)T_n(g)^*=A_1+A_2{\bf j}$ where 
$$
A_1={\bf I}_{n}-T_n(s)T_n(s)^*-T_n(h)T_n(h)^*, \quad A_2=T_n(s)T_n(h)^\top-T_n(h)T_n(s)^\top
$$
and the statement follows immediately, by Remark \ref{R:2.2}.\qed

\smallskip

Note that positive semidefiniteness of diagonal blocks in \eqref{3.5} is equivalent
to the inequality $|s(\zeta)|^2+|s(\zeta)|^2\le 1$ holding for all $\zeta\in\D$
which is necessary (since $|g(\zeta)|^2=|s(\zeta)|^2+|s(\zeta)|^2$) for $g$ to be in 
$\mathcal {QS}$ but not sufficient.

\bibliographystyle{amsplain}

\end{document}